\documentclass{amsproc}
\usepackage{graphicx}
\usepackage{amssymb}
\usepackage{epstopdf}
\usepackage{amsmath,amsthm,amscd,amssymb}
\usepackage{latexsym}
\usepackage[colorlinks,citecolor=red,pagebackref,hypertexnames=false]{hyperref}
\usepackage{geometry}                
\numberwithin{equation}{section}

\theoremstyle{plain}
\newtheorem{theorem}{Theorem}[section]
\newtheorem{lemma}[theorem]{Lemma}

\newtheorem{conjecture}[theorem]{Conjecture}

\theoremstyle{definition}
\newtheorem{definition}[theorem]{Definition}

\newtheorem{case[theorem]}{Case}

\theoremstyle{remark}
\newtheorem{remark}[theorem]{Remark}

\numberwithin{equation}{section}

\begin{document}

\title{On the unit distance problem} 

\author{A. Iosevich}

\date{today}

\email{iosevich@math.rochester.edu}

\address{Department of Mathematics, University of Rochester, Rochester, NY 14627}

\thanks{The author was partially supported by NSA H98230-15-1-0319}

\begin{abstract} The Erd\H os unit distance conjecture in the plane says that the number of pairs of points from a point set of size $n$ separated by a fixed (Euclidean) distance is $\leq C_{\epsilon} n^{1+\epsilon}$ for any $\epsilon>0$. The best known bound is $Cn^{\frac{4}{3}}$. We show that if the set under consideration is well-distributed and the fixed distance is much smaller than the diameter of the set, then the exponent $\frac{4}{3}$ is significantly improved. Corresponding results are also established in higher dimensions. The results are obtained by solving the corresponding continuous problem and using a continuous-to-discrete conversion mechanism.  The degree of sharpness of results is tested using the known results on the distribution of lattice points dilates of convex domains. 

We also introduce the following variant of the Erd\H os unit distance problem: how many pairs of points from a set of size $n$ are separated by an integer distance? We obtain some results in this direction and formulate a conjecture. 
\end{abstract} 

\maketitle


\section{Introduction}

\vskip.125in 

One of the hardest longstanding conjectures in extremal combinatorics is the Erd\H os unit distance conjecture (\cite{Erd45}, see also \cite{BMP05}). It says that if $P$ is a planar point set with $n$ points, then the number of pairs of elements of $P$ a fixed Euclidean distance apart is bounded by $C_{\epsilon}n^{1+\epsilon}$ for every $\epsilon>0$. The best known bound, obtain by Spencer, Szemeredi and Trotter (\cite{SST84}) is $Cn^{\frac{4}{3}}$. An interesting development occurred in 2005 when Pavel Valtr (\cite{V05}) proved that if the Euclidean distance is replaced by a distance induced by the norm defined by a bounded convex set with a smooth boundary and non-vanishing curvature, then the $Cn^{\frac{4}{3}}$ bound is, in general, best possible. 

The purpose of this paper is to show in the realm of well-distributed sets that the $Cn^{\frac{4}{3}}$ can be significantly improved if we count the number of pairs of points separated by a distance that is much smaller than the diameter of the set. Our main result is the following. 

\begin{definition} We say that $P \subset {\Bbb R}^d$ of size $n$ is {\it well-distributed} if there exists $c>0$ such that $|p-p'| \ge c$ and every unit lattice cube in ${\Bbb R}^d \cap {[0, n^{\frac{1}{d}}]}^d$ contains exactly one point of $P$. \end{definition} 

\vskip.125in  

\begin{theorem} \label{main} Let $B$ be a symmetric bounded convex set in ${\Bbb R}^d$, $d \ge 2$, with a smooth boundary and everywhere non-vanishing Gaussian curvature. Let $P$ be a well-distributed set of size $n$. Then for $k \in (1, n^{\frac{1}{d}})$, 
\begin{equation} \label{mamaest} \# \left\{(p,p') \in P \times P: k \leq {||p-p'||}_B \leq k+n^{-\frac{d-1}{d(d+1)}} \right\} \leq Cn^{2-\frac{2}{d+1}} \cdot \Lambda, \end{equation} where 
$$ \Lambda={\left( \frac{k}{n^{\frac{1}{d}}} \right)}^{\frac{d-1}{2}}.$$ \end{theorem} 

In particular, if $d=2$, the left hand side of (\ref{mamaest}) is bounded by $Cn^{\frac{4}{3}} \cdot {\left( \frac{k}{n^{\frac{1}{2}}} \right)}^{\frac{1}{2}}$, which is an improvement over the known $Cn^{\frac{4}{3}}$ bound when $k=o(n^{\frac{1}{2}})$. 

\vskip.125in 

\begin{remark} When $k \approx n^{\frac{1}{d}}$, Theorem \ref{main} is implicit in the main result in \cite{IJL09}, but the key feature here is the dependence on $k$ with the resulting improvement when $k=o(n^{\frac{1}{d}})$. Also, we shall prove below that in the case $k \approx n^{\frac{1}{d}}$, the estimate provided by Theorem \ref{main} is sharp. See also \cite{IS16} where the continuous-discrete correspondence is used in reverse in order to obtain sharpness examples for Falconer type estimates. \end{remark} 

\begin{remark} Note that the left hand side of (\ref{mamaest}) is trivially bounded by $Cn \cdot k^{d-1}$. Therefore, the estimate in Theorem \ref{main} is only interesting when $k>>n^{\frac{1}{d-1}-\frac{4}{(d-1)(d+1)}+\frac{1}{d(d-1)}}$. For example, in dimension two this threshold is $n^{\frac{1}{6}}$. \end{remark} 

\vskip.125in 

We also study the following variant of the Erd\H os unit distance problem. How many pairs of points from a set of $n$ points in ${\Bbb R}^d$, $d \ge 2$, are separated by an integer distance? When $P={\Bbb Z}^2 \cap {[0, n^{\frac{1}{d}}]}^d$, it is not difficult to see that the number of such pairs is 
$\approx n^{2-\frac{1}{d}}$. 

\begin{conjecture} \label{integerconjecture} Let $P \subset {\Bbb R}^d$, $d=2,3$, be a finite point set of size $n$. Then 
$$ \# \{(p,p') \in P \times P: |p-p'| \in {\Bbb Z} \} \leq Cn^{2-\frac{1}{d}}.$$ 
\end{conjecture} 

In higher dimensions this conjecture is not true, in general, due to the existence of the celebrated Lens example (see e.g. \cite{BMP05}) which shows that in dimensions $4$ and higher there exists $P \subset {\Bbb R}^d$ of size $n$ such that $\# \{(p,p') \in P \times P: |p-p'|=1 \} \ge cn^2$. But in the setting of well-distributed sets, Conjecture \ref{integerconjecture} still makes sense when $d \ge 4$. The following result follows easily from Theorem \ref{main}. 

\begin{theorem} \label{maininteger} Let $P$ be a well-distributed set of size $n$. Then 
$$ \# \{(p,p') \in P \times P: dist(|p-p'|, {\Bbb Z})<n^{-\frac{d-1}{d(d+1)}} \} \leq Cn^{2-\frac{2}{d+1}+\frac{1}{d}}=n^{2-\frac{1}{d}} \cdot n^{\frac{2}{d(d+1)}}.$$ 
\end{theorem} 

\vskip.25in 

\subsection{Sharpness of results} The results associated with the lattice point counting problems provide a useful tool for testing sharpness of Theorem \ref{main}. Let $n \approx q^d$ and $P={\Bbb Z}^d \cap B(\vec{0},10q)$, the ball of radius $10q$ centered at the origin. Let $N_d(R)$ denote the number of elements of ${\Bbb Z}^d$ inside the ball of radius $R$ centered at the origin. It is known (see e.g \cite{Hux96}) that 
$$N_d(R)=\omega_d R^d+D_d(R),$$ where $\omega_d$  is the volume of the unit ball, $|D_2(R)| \leq C_{\epsilon} R^{\frac{131}{208}+\epsilon}$ (\cite{Hux03}), $|D_3(R)| \leq C_{\epsilon} R^{\frac{21}{16}+\epsilon}$ (\cite{HB99}), and $|D_d(R)| \leq C_{\epsilon} R^{d-2+\epsilon}$ for $d \ge 4$ (\cite{Fr82}).

Then 
$$ \# \left\{(p,p') \in P \times P: q \leq |p-p'| \leq q+q^{-\frac{d-1}{d+1}} \right\} \ge Cq^d \cdot \left(N\left(q+q^{-\frac{d-1}{d+1}}\right)-N(q)\right).$$ 

We have 
$$ N\left(q+q^{-\frac{d-1}{d+1}} \right)-N(q)=\omega_d \left({\left(q+q^{-\frac{d-1}{d+1}}\right)}^d-q^d \right)+D \left(q+q^{-\frac{d-1}{d+1}}\right)-D(q).$$ 

Using the bounds on $|D(R)|$ described above, we see that 
$$ N\left(q+q^{-\frac{d-1}{d+1}}\right)-N(q) \ge cq^{d-1-\frac{d-1}{d+1}}, $$ which implies that 
$$  \# \{(p,p') \in P \times P: q \leq |p-p'| \leq q+q^{-\frac{d-1}{d+1}} \} \ge Cq^{2d-1+\frac{d-1}{d+1}} \approx n^{2-\frac{2}{d+1}},$$ proving that Theorem \ref{main} is sharp when $k \approx n^{\frac{1}{d}}$. 

When $k<<q$, we can conclude that 
$$ N\left(k+q^{-\frac{d-1}{d+1}}\right)-N(k) \ge ck^{d-1} q^{-\frac{d-1}{d+1}}$$ if the right hand side is larger than the error term measured in terms of the bounds on $|D_d(R)|$ described above. This happens for a range of $k$'s. In this range, 
$$ \# \{(p,p') \in P \times P: q \leq |p-p'| \leq k+q^{-\frac{d-1}{d+1}} \} \ge Ck^{d-1} q^d q^{-\frac{d-1}{d+1}}.$$ 

The right hand side is smaller then the bound obtained by Theorem \ref{main} when $k=o(n^{\frac{1}{d}})$. This may indicate that Theorem \ref{main} is not sharp in this range, but it is also possible that a more sophisticated sharpness example may be found. 

\vskip.125in 

The construction above applied to $\approx n^{\frac{1}{d}}$ annuli shows that when $k \approx n^{\frac{1}{d}}$, the conclusion of Theorem \ref{maininteger} is sharp. 

\vskip.125in

\subsection{Acknowledgements}  The author wishes to thank Adam Sheffer and Josh Zahl for some very helpful remarks and suggestions. 

\vskip.125in 

\section{Proof of the main result} 

\vskip.125in 

For $\frac{d}{2}<s<d$, define
$$ \mu_{q,s}(x)=q^{-d} q^{\frac{d^2}{s}} \sum_{p \in P} \phi \left(q^{\frac{d}{s}}\left(x-\frac{p}{q}\right)\right) \phi\left(\frac{p}{q}\right),$$ where $\phi$ is a smooth cut-off function supported in the ball of radius $2$ and identically equal to $1$ in the ball of radius $1$. 

This is a natural measure on the $q^{-\frac{d}{s}}$-neighborhood of $\frac{1}{q}P=\left\{ \frac{a}{q}: a \in P \right\}$. Our goal is to bound the expression 
\begin{equation} \label{keyexpression} \int \int_{\left\{(x,y): t \leq {||x-y||}_B \leq t+q^{-\frac{d}{s}}\right\}} d\mu_{q,s}(x) d\mu_{q,s}(y), \end{equation} where ${||\cdot||}_B$ is the norm induced by a bounded symmetric convex set $B$ with a smooth boundary and everywhere non-vanishing curvature, and then relate it to the count for the number of pairs separated by a given distance. 

Using a Fourier inversion type argument (see e.g. \cite{M16} or \cite{W03}), the expression (\ref{keyexpression}) equals a constant multiple of 
\begin{equation} \label{keyfourierexpression} \int {|\widehat{\mu}_{q,s}(\xi)|}^2 \widehat{\chi}_{A_{t,q,s}}(\xi) d\xi, \end{equation} where 
$$ A_{t,q,s}=\{x \in {\Bbb R}^d: t \leq {||x||}_B \leq t+q^{-\frac{d}{s}}\},$$ $\chi$ denotes its indicator function, ${||\cdot||}_B$ is the norm induced by a symmetric bounded convex set $B$ with a smooth boundary and everywhere non-vanishing curvature and 
$$ \widehat{f}(\xi)=\int e^{-2 \pi i x \cdot \xi} f(x)dx,$$ defined for functions in $L^2({\Bbb R}^d)$. 

\vskip.125in 

We first show that the $s$-energy integral of $\mu_{q,s}$ is bounded independently of $q$. See, for example, \cite{Falc86II}, \cite{M95} and \cite{M16} for the background on energy integrals in the setting of sets of a given Hausdorff dimension. 

\begin{lemma} \label{energylemma} For any $s \in \left(\frac{d}{2},d \right)$, 
$$ \int {|\widehat{\mu}_{q,s}(\xi)|}^2 {|\xi|}^{-d+s} d\xi=c_{d,s} 
\int \int {|x-y|}^{-s} d\mu_{q,s}(x) d\mu_{q,s}(y) \leq C<\infty$$ with a bound independent of $q$. 
\end{lemma} 

We shall give the proof of Lemma \ref{energylemma} at the end of the paper. Next, we bound the Fourier transform of the indicator function of $A_{t,q,s}$. 

\begin{lemma} \label{fourierannuluslemma} (\cite{Falc86}) With the notation above, 
\begin{equation} \label{fourierannulusest} \left| \widehat{\chi}_{A_{t,q,s}}(\xi) \right| \leq Ct^{\frac{d-1}{2}} {|\xi|}^{-\frac{d-1}{2}} \min \left\{ q^{-\frac{d}{s}}, {|\xi|}^{-1} \right\}, \end{equation} where $C$ is a universal constant independent of $t$ or $q$. 
\end{lemma} 

Falconer proved this result in (\cite{Falc86}) in the case when $B$ is the unit ball. The proof of the general case is similar. 

\vskip.125in 

With Lemma \ref{energylemma} and Lemma \ref{fourierannuluslemma} in tow, we see that 
\begin{equation} \label{keyprepest} \int {|\widehat{\mu}_{q,s}(\xi)|}^2 \widehat{\chi}_{A_{t,q,s}}(\xi) d\xi \leq Ct^{\frac{d-1}{2}} \cdot q^{-\frac{d}{s}} 
\int {|\widehat{\mu}_{q,s}(\xi)|}^2 {|\xi|}^{-\frac{d-1}{2}} d\xi \end{equation}
$$ \leq C't^{\frac{d-1}{2}} q^{-\frac{d}{s}} \int {|\widehat{\mu}_{q,s}(\xi)|}^2 {|\xi|}^{-d+s} d\xi \leq C''t^{\frac{d-1}{2}} q^{-\frac{d}{s}}$$ if 
$s \ge \frac{d+1}{2}$.

\vskip.125in 

We are now ready for the combinatorial conclusion. See \cite{HI05}, \cite{IL05}, \cite{I08} and \cite{IRU14} where various forms of the continuous to discrete conversion mechanisms are developed and applied. Observe that 
$$ \# \{(p,p') \in P \times P: k \leq {||x-y||}_B \leq k+q^{-\frac{d}{s}+1} \}$$
$$ \leq C q^{2d} \cdot \int \int_{\left\{(x,y): \frac{k}{q} \leq {||x-y||}_B \leq \frac{k}{q}+q^{-\frac{d}{s}}\right\}} d\mu_{q,s}(x) d\mu_{q,s}(y).$$

By (\ref{keyprepest}) this expression is bounded by 
$$ C' q^{2d-\frac{d}{s}} q^{-\frac{d-1}{2}} k^{\frac{d-1}{2}}=C' n^{2-\frac{2}{d+1}} \cdot {\left( \frac{k}{n^{\frac{1}{d}}} \right)}^{\frac{d-1}{2}},$$ as desired. This completes the proof of Theorem \ref{main} up the proof of Lemma \ref{energylemma}. 

\vskip.125in 

\subsection{Proof of Lemma \ref{energylemma}} This result is proved in \cite{IRU14}, but we include the proof for the sake of completeness. We have 
$$ q^{-2d} q^{\frac{2d^2}{s}} \int \int {|x-y|}^{-s} d\mu_{q,s}(x) d\mu_{q,s}(y)$$
$$=\sum_{p,p' \in P} \phi\left(\frac{p}{q}\right) \phi\left(\frac{p'}{q}\right) 
\int \int {|x-y|}^{-s} \phi \left(q^{\frac{d}{s}}\left(x-\frac{p}{q}\right)\right) \phi \left(q^{\frac{d}{s}}\left(y-\frac{p'}{q}\right)\right) dxdy=I+II$$ where 
$$ I=q^{-2d} q^{\frac{2d^2}{s}} \sum_{p \in P} \phi^2\left(\frac{p}{q}\right) \int \int {|x-y|}^{-s} \phi \left(q^{\frac{d}{s}}\left(x-\frac{p}{q}\right)\right) \phi \left(q^{\frac{d}{s}}\left(y-\frac{p}{q}\right)\right) dxdy$$ and
$$II=q^{-2d} q^{\frac{2d^2}{s}} \sum_{p \not=p' \in P} \phi\left(\frac{p}{q}\right) \phi\left(\frac{p'}{q}\right) 
\int \int {|x-y|}^{-s} \phi \left(q^{\frac{d}{s}}\left(x-\frac{p}{q}\right)\right) \phi \left(q^{\frac{d}{s}}\left(y-\frac{p'}{q}\right)\right) dxdy.$$

By a direct calculation, $I$ is bounded. Using the separation between $p$ and $p'$, we see that 
$$ II \leq q^{-2d} q^{\frac{2d^2}{s}} \sum_{p \not=p' \in P} \phi\left(\frac{p}{q}\right) \phi\left(\frac{p'}{q}\right) {\left| \frac{p}{q}-\frac{p'}{q} \right|}^{-s} \cdot q^{\frac{-2d^2}{s}}$$
$$=q^{-2d} q^s \sum_{p \not=p' \in P} \phi\left(\frac{p}{q}\right) \phi\left(\frac{p'}{q}\right) {|p-p'|}^{-s}.$$

Using the well-distributivity assumption on $P$ we may replace the sum by the integral and thus the quantity is bounded. This completes the proof of Lemma \ref{energylemma}. 

\vskip.25in 

\section{Proof of Theorem \ref{maininteger}} 

\vskip.125in 

We have shown above that 
$$ \# \{(p,p') \in P \times P: k \leq |p-p'| \leq k+n^{-\frac{d-1}{d(d+1)}} \} \leq C n^{2-\frac{2}{d+1}} \cdot {\left( \frac{k}{n^{\frac{1}{d}}} \right)}^{\frac{d-1}{2}}.$$ 

Summing both sides over $k=1,2 \dots, N$, where $N \approx n^{\frac{1}{d}}$ yields the conclusion of Theorem \ref{maininteger}.

\bigskip


\begin{thebibliography}{8}




\bibitem{BMP05} P. Brass, W. O. Moser, J. Pach. {\em Research Problems in Discrete Geometry.} Springer, (2005).


\bibitem{Erd45} P. Erd\H os. {\it On sets of distances of n points} Amer. Math. Monthly. \textbf{53} (1946), 248--250.




\bibitem{Falc86} K. J. Falconer, {\it On the Hausdorff dimensions of distance sets}, Mathematika \textbf{32} (1986), 206-212.

\bibitem{Falc86II} K. J. Falconer, {\it The geometry of fractal sets}, Cambridge Tracts in Mathematics \textbf{85}, Cambridge Univ. Press, Cambridge, (1986). 

\bibitem{Fr82} F. Fricker, {\it Einfuhrung die Gitterpunktlehre}, Birkhauser, Verlag, (1982). 






\bibitem{HB99} D. R. Heath-Brown, {\it Lattice points in the sphere}, Number theory in progress, Vol. 2 (Zakopane-Ko?cielisko, 1997), 883-892, de Gruyter, Berlin, (1999).

\bibitem{HI05} S. Hofmann and A. Iosevich {\it Circular averages and Falconer/Erd\"os distance conjecture in the plane for random metrics} Proc. Amer. Mat. Soc. \textbf{133} (2005) 133-144. 

\bibitem{Hux96} M. N. Huxley, {\it Area, Lattice Points, and Exponential Sums}, London Mathematical Society Monographs New Series \textbf{13}, Oxford Univ. Press, (1996).

\bibitem{Hux03} M. N. Huxley, {\it Exponential sums and lattice points. III} Proc. London Math. Soc. (3) 87 (2003), no. 3, 591-609.


\bibitem{I08} A. Iosevich, {\it Fourier analysis and geometric combinatorics.} Topics in Mathematical Analysis, Series on Analysis, Application and Computation, \textbf{3}, World Scientific, proceedings of the Padova lectures in analysis in 2004 and 2005 (2008). 

\bibitem{IJL09} A. Iosevich, H. Jorati and I. Laba, {\it Geometric incidence theorems via Fourier analysis}, Trans. Amer. Math. Soc. \textbf{361} (2009), no. 12, 6595-6611.

\bibitem{IL05} A. Iosevich, I. {\L}aba, {\sl K-distance sets, Falconer conjecture, and discrete analogs},

\bibitem{IRU14} A. Iosevich, M. Rudnev and I. Uriarte-Tuero, {\it Theory of dimension for large discrete sets and applications}, Math. Model. Nat. Phenom. \textbf{9} (2014), no. 5, 148-169.

\bibitem{IS16} A. Iosevich and S. Senger, {\it Sharpness of Falconer's $\frac{d+1}{2}$ estimate}, Ann. Acad. Sci. Fenn. Math. \textbf{41} (2016), no. 2, 713-720. 











\bibitem{M95} P. Mattila, {\it Geometry of sets and measures in Euclidean spaces}, Cambridge University Press, (1995). 

\bibitem{M16} P. Mattila, {\it Fourier Analysis and Hausdorff dimension}, Cambridge University Press, Cambridge studies in advanced mathematics, \textbf{150}, (2016). 













\bibitem{SST84} J. Spencer, E. Szemer\'edi, and W. T. Trotter. {\it Unit distances in the Euclidean plane} B. Bollob\'as, editor, ``Graph Theory and Combinatorics," pages 293-303. Academic Press, New York, NY, 1984.


\bibitem{V05} P. Valtr, {\it Strictly convex norms allowing many unit distances and related touching questions}, manuscript (2005).



\bibitem{W03} T. Wolff, {\it Lectures on harmonic analysis}, I. Laba and C. Shubin, eds. University Lecture Series, \textbf{29}. Amer. Math. Soc., Providence, RI, (2003).




\end{thebibliography}
\end{document}